\documentclass[12pt]{article}
\usepackage{amsfonts,amssymb}
\newcommand{\cartht}[2]{\left[\begin{array}{c}{#1}\\{#2}\end{array}\right]}
\newcommand{\sfrac}[2]{{\textstyle{\frac{#1}{#2}}}}

\begin{document}
\newpage
\pagestyle{empty}
\setcounter{page}{0}
\vfill
\begin{center}

{\Large \textbf{$q$-deformed ${\cal W}$-algebras and elliptic algebras}}

\vspace{10mm}

{\large J. Avan}

\vspace{4mm}

{\em LPTHE, CNRS-URA 280, Universit{\'e}s Paris VI/VII, France}

\vspace{7mm}

{\large L. Frappat, M. Rossi, P. Sorba}

\vspace{4mm}

{\em Laboratoire d'Annecy-le-Vieux de Physique Th{\'e}orique LAPTH, 
CNRS-URA 1436}

{\em LAPP, BP 110, F-74941 Annecy-le-Vieux Cedex, France}

\end{center}

\vfill
\vfill

\begin{abstract}
The elliptic algebra ${\cal A}_{q,p}(\widehat{sl}(N)_{c})$ at the 
critical level $c=-N$ has an extended center containing trace-like 
operators $t(z)$.  Families of Poisson structures, defining 
$q$-deformations of the ${\cal W}_{N}$ algebra, are constructed.  The 
operators $t(z)$ also close an exchange algebra when $(-p^{1/2})^{NM} = 
q^{-c-N}$ for $M \in {\mathbb Z}$.  It becomes Abelian when in addition 
$p=q^{Nh}$ where $h$ is a non-zero integer.  The Poisson structures 
obtained in these classical limits contain different $q$-deformed ${\cal 
W}_{N}$ algebras depending on the parity of $h$, characterizing the 
exchange structures at $p \ne q^{Nh}$ as new ${\cal W}_{q,p}(sl(N))$ 
algebras.  
\end{abstract}

\vfill
\vfill

\vfill
\emph{Talk presented by L. Frappat at the 7th International Colloquium 
``Quantum Groups and Integrable Systems'' (QGIS), Prague, 18-20 June 1998.}
\vfill

\rightline{math.QA/9807131}
\rightline{July 1998}

\newpage
\pagestyle{plain}

\section{Introduction}

The results hereafter presented are based on four papers of the authors 
[1--3].  They deal with the construction of $q$-deformed Virasoro, and 
more generally ${\cal W}$, algebras in relation with the quantum 
elliptic algebras ${\cal A}_{q,p}(\widehat{sl}(N)_{c})$.  Indeed what we 
achieve here is really a $q$-deformation of the Sugawara construction 
\emph{on an abstract level} (that is only abstract algebraic relations 
in ${\cal A}_{q,p}(\widehat{sl}(N)_{c})$ are used).

We first prove the existence of an extended center in ${\cal 
A}_{q,p}(\widehat{sl}(N)_{c})$ at $c=-N$ and compute the Poisson 
structures on this center.  These structures are identified as 
$q$-deformed $Vir_{q}(sl(N))$ algebras.  We then show the existence of 
closed (quadratic) exchange algebras whenever $(-p^\frac{1}{2})^{NM} = 
q^{-c-N}$ for any integer $M \in {\mathbb Z}$.  These algebras admit a 
classical limit (commuting algebras) at $p = q^{Nh}$ with $h \in 
{\mathbb Z} \backslash \{ 0 \}$.  The related Poisson structures include 
for $h$ even the Poisson structures in \cite{FR}.  The exchange algebras 
therefore realize new quantizations of these Poisson structures.  When 
$h$ is odd, by contrast, this classical limit takes a form different 
from the initial $Vir_{q}(sl(N))$ structures.  This emphasizes the key 
role of the initial 3-parameter structure ${\cal 
A}_{q,p}(\widehat{sl}(N)_{c})$ in allowing for an intermediate quantum 
$q$-deformed Virasoro algebra.  When computing the mode expansion of 
these structures (both classical and quantum), a ``sector-type'' 
labeling is needed due to the singularity behaviour of the structure 
function viewed as an analytic continuation and therefore exhibiting 
different formal series expansions corresponding to different convergent 
series expansions in distinct domains of the complex plane.  Finally, we 
describe the essential features of the extension of our computation to 
higher spin fields.

\section{The elliptic algebra ${\cal A}_{q,p}(\widehat{sl}(N)_{c})$}

We start by defining the $R$-matrix of the ${\mathbb Z}_{N}$-vertex model 
(${\mathbb Z}_N \equiv {\mathbb Z}/N{\mathbb Z}$) \cite{Bax,Bela}:
\[
\widetilde R(z,q,p) = \frac{z^{2/N-2}}{\kappa(z^2)} 
\frac{\vartheta\cartht{\textstyle{\frac{1}{2}}} 
{\textstyle{\frac{1}{2}}}(\zeta,\tau)} 
{\vartheta\cartht{\textstyle{\frac{1}{2}}} 
{\textstyle{\frac{1}{2}}}(\xi+\zeta,\tau)} \,\, 
\sum_{(\alpha_1,\alpha_2)} W_{(\alpha_1,\alpha_2)}(\xi,\zeta,\tau) \,\, 
I_{(\alpha_1,\alpha_2)} \otimes I_{(\alpha_1,\alpha_2)}^{-1} \,,
\]
where $z=e^{i\pi\xi}$, $q=e^{i\pi\zeta}$, $p=e^{2i\pi\tau}$ and 
$(\alpha_1,\alpha_2)\in{\mathbb Z}_N\times{\mathbb Z}_N$.  \\
$\vartheta$ are the Jacobi theta functions with rational characteristics 
$(\gamma_1,\gamma_2) \in {\mathbb Z}/N \times {\mathbb Z}/N$:
\[
\vartheta\cartht{\gamma_1}{\gamma_2}(\xi,\tau) = \sum_{m \in {\mathbb Z}}
\exp\Big(i\pi(m+\gamma_1)^2\tau + 2i\pi(m+\gamma_1)(\xi+\gamma_2)
\Big) \,.
\]
The normalization factor is given by:
\[
\frac{1}{\kappa(z^2)} = \frac{(q^{2N}z^{-2};p,q^{2N})_\infty
\, (q^2z^2;p,q^{2N})_\infty \, (pz^{-2};p,q^{2N})_\infty \,
(pq^{2N-2}z^2;p,q^{2N})_\infty} {(q^{2N}z^2;p,q^{2N})_\infty
\, (q^2z^{-2};p,q^{2N})_\infty \, (pz^2;p,q^{2N})_\infty \,
(pq^{2N-2}z^{-2};p,q^{2N})_\infty} \,.
\]
The functions $W_{(\alpha_1,\alpha_2)}$ are defined as follows:
\[
W_{(\alpha_1,\alpha_2)}(\xi,\zeta,\tau) = \frac{1}{N} \,\, \frac
{\vartheta\cartht{\textstyle{\frac{1}{2}}+\alpha_1/N}
{\textstyle{\frac{1}{2}}+\alpha_2/N}(\xi+\zeta/N,\tau)} 
{\vartheta\cartht{\textstyle{\frac{1}{2}}+\alpha_1/N}
{\textstyle{\frac{1}{2}}+\alpha_2/N}(\zeta/N,\tau)} \,,
\]
and the matrices $I_{(\alpha_1,\alpha_2)}$ by $I_{(\alpha_1,\alpha_2)} = 
g^{\alpha_{2}} \, h^{\alpha_{1}}$ where $g_{ij} = 
\omega^{i}\delta_{ij}$, $h_{ij} = \delta_{i+1,j}$ are $N \times N$ 
matrices (the addition of indices being understood modulo $N$) and 
$\omega = e^{2i\pi/N}$.  \\
The $R$-matrix $\widetilde R$ is ${\mathbb Z}_N$-symmetric:
\[
\widetilde R_{c+s\,,\,d+s}^{a+s\,,\,b+s} = \widetilde 
R_{c\,,\,d}^{a\,,\,b} \qquad a,b,c,d,s \in {\mathbb Z}_N \,.
\]
We introduce the ``gauge-transformed'' matrix:
\[
R(z,q,p) = (g^{\frac{1}{2}} \otimes g^{\frac{1}{2}}) \widetilde
R(z,q,p) (g^{-\frac{1}{2}} \otimes g^{-\frac{1}{2}}) \,.
\]
It satisfies the following properties: \\
-- Yang--Baxter equation: $R_{12}(z) \, R_{13}(w) \, R_{23}(w/z) = 
R_{23}(w/z) \, R_{13}(w) \, R_{12}(z)$, \\
-- unitarity: $R_{12}(z) R_{21}(z^{-1}) = 1$, \\
-- crossing symmetry: $R_{12}(z)^{t_2} R_{21}(q^{-N}z^{-1})^{t_2} = 1$, 
\\
-- antisymmetry: $R_{12}(-z) = \omega \, (g^{-1} \otimes {\mathbb I}) \, 
R_{12}(z) \, (g \otimes {\mathbb I})$, \\
-- quasi-periodicity: $\widehat{R}_{12}(-p^{\frac{1}{2}}z) = 
(g^{\frac{1}{2}}hg^{\frac{1}{2}} \otimes {\mathbb I})^{-1} \left 
(\widehat{R}_{21}(z^{-1})\right )^{-1} (g^{\frac{1}{2}}hg^{\frac{1}{2}} 
\otimes {\mathbb I})$, \\ 
where $\widehat{R}_{12}(x) = \tau_{N}(q^{1/2}x^{-1}) \, R_{12}(x)$ and 
$\tau_N(z) = z^{\frac{2}{N}-2} \, 
\frac{\Theta_{q^{2N}}(qz^2)}{\Theta_{q^{2N}}(qz^{-2})}$.  The function 
$\tau_N(z)$ is periodic with period $q^N$: $\tau_N(q^Nz) = \tau_N(z)$ 
and satisfies $\tau_N(z^{-1}) = \tau_N(z)^{-1}$.

\medskip

We now define the elliptic quantum algebra ${\cal 
A}_{q,p}(\widehat{sl}(N)_{c})$ \cite{FIJKMY,JKOS}.  It is an algebra of 
operators $L_{ij}(z) \equiv \sum_{n \in {\mathbb Z}} L_{ij}(n) \, z^{n}$ 
where $i,j \in {\mathbb Z}_{N}$:
\[
L(z) = \left(\begin{array}{ccc} L_{11}(z) & \cdots &
L_{1N}(z) \cr \vdots && \vdots \cr L_{N1}(z) & \cdots &
L_{NN}(z) \cr \end{array}\right) \,.
\]
The $q$-determinant is given by ($\varepsilon(\sigma)$ being the 
signature of the permutation $\sigma$):
\[
q\mbox{-}\det L(z) \equiv \displaystyle \sum_{\sigma\in{\mathfrak S}_N}
\varepsilon(\sigma) \prod_{i=1}^N L_{i,\sigma(i)}(z q^{i-N-1}) \,.
\]
${\cal A}_{q,p}(\widehat{gl}(N)_c)$ is defined by imposing the following 
constraints on $L(z)$:
\[
\widehat{R}_{12}(z/w) \, L_1(z) \, L_2(w) = L_2(w) \, L_1(z) \,
\widehat{R}_{12}^{*}(z/w) \,,
\label{eq218}
\]
where $L_1(z) \equiv L(z) \otimes {\mathbb I}$, $L_2(z) \equiv {\mathbb 
I} \otimes L(z)$ and $\widehat R^{*}_{12}(z,q,p) \equiv \widehat 
R_{12}(z,q,p^*=pq^{-2c})$.  \\
The $q$-determinant is in the center of ${\cal A}_{q,p}(\widehat{gl}(N)_c)$
and one sets
\[
{\cal A}_{q,p}(\widehat{sl}(N)_c) = {\cal A}_{q,p}(\widehat{gl}(N)_c)/
\langle q\mbox{-}\det L - q^{\frac{c}{2}} \rangle \,.
\]

\section{The center of ${\cal A}_{q,p}(\widehat{sl}(N)_{c})$}

We set $L^+(z) \equiv L(q^{c/2}z)$ and $L^-(z) \equiv 
g^{\frac{1}{2}}hg^{\frac{1}{2}} \, L(-p^{\frac{1}{2}}z) \, 
(g^{\frac{1}{2}}hg^{\frac{1}{2}})^{-1}$.  Let us define the following 
operator in ${\cal A}_{q,p}(\widehat{sl}(N)_{c})$:
\[
t(z) \equiv {\rm Tr}({\cal L}(z)) \equiv 
{\rm Tr} \left(L^+(q^{c/2}z) \, L^-(z)^{-1}\right) \,.
\]
{From} the defining relations of the elliptic algebra and the properties 
of the matrix $\widehat R$ (especially the quasi-periodicity 
property), one gets the following theorems:  \\
\textbf{Theorem 1}:  
At the critical level $c=-N$, the operators generated by $t(z)$ lie in 
the center of the algebra ${\cal A}_{q,p}(\widehat{sl}(N)_{c})$:
\[
\left[t(z),L^+(w)\right] = \left[t(z),L^-(w)\right] = 0 \,.
\]
\textbf{Theorem 2}:  
At the critical level $c=-N$, the operators generated by $t(z)$ close an 
Abelian algebra:
\[
[t(z),t(w)] = 0 \,.
\]
This allows us to define a Poisson structure on the center of ${\cal 
A}_{q,p}(\widehat{sl}(N)_{c})$ at $c=-N$.  When $c \ne -N$, one has:
\begin{eqnarray*}
t(z)t(w) &=& {\cal Y}(w/z)^{i_1i_2}_{j_1j_2} \,\, {\cal L}(w)^{j_2}_{i_2} 
\,\, {\cal L}(z)^{j_1}_{i_1} \\
&=& t(w)t(z) + (c+N) \left( \frac{d{\cal Y}}{dc} \left(\frac{w}{z}\right) 
\right)^{i_1i_2}_{j_1j_2} \,\, \left.  {\cal L}(w)^{j_2}_{i_2} \,\, 
{\cal L}(z)^{j_1}_{i_1} \right\vert_{cr} + \cdots
\end{eqnarray*}
The Poisson structure is then defined by:
\[
\{ t(z),t(w) \} = \left( \frac{d{\cal Y}}{dc} \left(\frac{w}{z}\right) 
\right)^{i_1i_2}_{j_1j_2} \, \left.  {\cal L}(w)^{j_2}_{i_2} \, 
{\cal L}(z)^{j_1}_{i_1} \right\vert_{cr} \,.
\]
Decomposing ${\cal Y}(x)$ as ${\cal Y}(x) = T(x) \, {\cal M}(x)$ where
\[
T(x) = \frac{\tau_N(q^{\frac{1}{2}}x^{-1})\tau_N(q^{\frac{1}{2} - c}x)}
{\tau_N(q^{\frac{1}{2}}x)\tau_N(q^{\frac{1}{2} - c}x^{-1})}
\]
and
\[
{\cal M}(x) = \left( \left( R_{21}(x) \,
{R_{21}(q^{c+N}x)}^{-1} \, {R_{12}(x^{-1})}^{-1} \right)^{t_2}
\, {R_{12}(q^cx^{-1})}^{t_2} \right)^{t_2} \,, 
\]
one gets
\[
\left.  \frac{d{\cal Y}}{dc}(x) \right\vert_{cr} = 
\left.  \frac{dT}{dc}(x) \right\vert_{cr} {\mathbb I}_2 \otimes {\mathbb I}_2 \,\, 
+ \,\, \left.  \frac{d{\cal M}}{dc}(x) \right\vert_{cr} \,,
\]
since ${T(x)}_{cr} = 1$ and ${{\cal M}(x)}_{cr} = {\mathbb I}_2 \otimes 
{\mathbb I}_2$ at the critical level $c=-N$.  \\
The properties of the matrix $R_{12}$ imply the vanishing of the 
derivative of ${\cal M}$ with respect to $c$ at the critical level and 
it follows that:
\[
\{ t(z),t(w) \} = \left.  \frac{dT}{dc} \left(\frac{w}{z}\right) 
\right\vert_{cr} t(z) \, t(w) \equiv f(w/z) \, t(z) \, t(w) 
\]
where
\begin{eqnarray*}
f(x) &=& -2\ln q \left[ \sum_{l \ge 0}\, \left( 
\frac{2x^2q^{2Nl}}{1-x^2q^{2Nl}} - 
\frac{x^2q^{2Nl+2}}{1-x^2q^{2Nl+2}} -
\frac{x^2q^{2Nl-2}}{1-x^2q^{2Nl-2}} \right) \right. \\
&& \hspace{12mm} \left.  - \frac{x^2}{1-x^2} + {\textstyle{\frac{1}{2}}} 
\, \frac{x^2q^2}{1-x^2q^2} + {\textstyle{\frac{1}{2}}} \, 
\frac{x^2q^{-2}}{1-x^2q^{-2}} - (x \leftrightarrow x^{-1}) \right] \,.
\end{eqnarray*}
The Poisson structure on the center of ${\cal 
A}_{q,p}(\widehat{sl}(N)_{c})$ at $c=-N$ \emph{is not $p$ depending}.

\medskip

We now define the Poisson structure for modes of the generating function 
$f(x)$.  The modes of $t(z)$ are defined in the sense of formal series 
expansions: $\displaystyle t_m = \oint_{C} \frac{dz}{2i\pi z} \, z^{-m} 
\, t(z)$, where $C$ is a contour encircling the origin.  \\
However, the function $f(x)$ has poles for $x = \pm q^{\pm Nl}$ and $x = 
\pm q^{\pm Nl \pm 1}$ with $l \ge 0$, therefore
\[
\oint_{C_1} \frac{dz}{2\pi iz} \oint_{C_2} 
\frac{dw}{2\pi iw} \ne \oint_{C_2} \frac{dz}{2\pi iz} \oint_{C_1} 
\frac{dw}{2\pi iw} \,.
\]
The relative position of the contours $C_{1}$ and $C_{2}$ must be 
specified in order to have an unambiguous result for the Poisson 
bracket.  Moreover, the antisymmetry of the Poisson bracket is only 
guaranteed at the mode level by an explicit symmetrization with respect 
to the position of the contours $C_{1}$ and $C_{2}$.  The mode Poisson 
bracket is thus defined as:
\begin{eqnarray*}
\{ t_n,t_m \} &=& {\textstyle{\frac{1}{2}}} \left( \oint_{C_1} 
\frac{dz}{2\pi iz} \oint_{C_2} \frac{dw}{2\pi iw} + \oint_{C_2} 
\frac{dz}{2\pi iz} \oint_{C_1} \frac{dw}{2\pi iw} \right) \\
&& \hspace{30mm} z^{-n}w^{-m} \, f(w/z) \, t(z) \, t(w)
\end{eqnarray*}
where $C_{1}$ and $C_{2}$ are circles of radii $R_{1}$ and $R_{2}$ and 
one chooses $R_{1} > R_{2}$.

\medskip

The Poisson bracket depends on the relative positions of the contours 
$C_1$ and $C_2$, circles of radii $R_1$ and $R_2$.  Thus we introduce 
the notion of sectors: the sector $(k)$ is defined by $R_{1}/R_{2} \in 
\left] {\vert q\vert }^{-P(k)},{\vert q\vert }^{-P(k+1)} \right[ \,$ 
where $\{ P(k) \}$ is the ordered set of powers of $q^{-1}$ where the 
poles of $f$ are located: $\{ P(k) \} = \{0,1; N-1,N,N+1; \dots; 
Nl-1,Nl,Nl+1; \dots\}$.

In the case of $sl(N)$, for $R_1/R_2 \in ]1,\vert q \vert^{-1}[$ (first 
sector), one gets the following formula:
\[
\{ t_n,t_m \}_{k=0} = -2(q-q^{-1}) \ln q \sum_{r\in{\mathbb Z}}
\frac{[(N-1)r]_{q}[r]_{q}}{[Nr]_{q}} \, t_{n-2r} \, t_{m+2r} 
\]
where $[r]_{q}$ are $q$-numbers: $\displaystyle [r]_{q} \equiv 
\frac{q^r-q^{-r}}{q-q^{-1}}$.

In the case of $sl(2)$, for $R_1/R_2 \in ]\vert q \vert^{\pm k},\vert q 
\vert^{\pm (k+1)}[$, it is possible to get a nice compact formula for 
any sector:
\begin{eqnarray*}
\{ t_n,t_m \}_k &=& (-1)^{k+1} 2 \ln q \oint_{C_1} \frac{dz}{2\pi iz} 
\oint_{C_2} \frac{dw}{2\pi iw} \\
&& \sum_{s \in {\mathbb Z}} \frac{q^{(2k+1)s} - q^{-(2k+1)s}}{q^s+q^{-s}} \,
(\frac zw)^{2s} \, z^{-n}w^{-m} \, t(z) \, t(w) \,.
\end{eqnarray*}

It follows that one obtains in this way a family of Poisson structures 
indexed by $k \in {\mathbb Z}$, the sector $k=0$ corresponding to the 
Poisson bracket in \cite{FR}.

\medskip

To realize deformed ${\cal W}_{N}$ Poisson structures, we need to 
introduce generating functions for the higher spin objects.  Having at 
our disposal only one commuting generating function $t(z)$, we are led 
to define shifted products, although with the same generator.  We define 
therefore the higher spin objects for $1 \le i \le N-1$ as:
\[
s_{i}(z) = \prod_{u=-(i-1)/2}^{(i-1)/2} t(q^uz)
\]
The Poisson algebra reads as:
\[
\{ s_{i}(z),s_{j}(w) \} = \sum_{u=-(i-1)/2}^{(i-1)/2} 
\sum_{v=-(j-1)/2}^{(j-1)/2} f\Big(q^{v-u}\frac{w}{z}\Big) \, s_{i}(z) \, 
s_{j}(w) \,.
\]
As before, one can introduce the notion of sectors for the Poisson 
brackets of the mode expansions, now labelled by $k(i,j)$ \cite{AFRS3}.  \\
{For} the sector $k = 0$, one gets:
\begin{eqnarray*}
\{ s_i(n),s_j(m) \} &=& - \sum_{r\in{\mathbb Z}} 
\frac{[(N-\max(i,j))r]_{q}[\min(i,j)r]_{q}}{[Nr]_{q}} \nonumber \\
&&\times 2(q-q^{-1}) \ln q \, s_i(n-2r) s_j(m+2r) \,. \qquad \qquad (*)
\end{eqnarray*}
Any Poisson structure in a given sector $k(i,j)$ can be obtained from 
(*) by adding to the $r$-dependent structure coefficient 
contributions from the relevant singularities of the structure function, 
i.e. power series expansions of terms $\delta(\pm q^{\pm P_{ij}(s)}w/z) 
\, s_{i}(z)$ $s_{j}(w)$ for $1 \le s \le k(i,j)$.

\section{${\cal W}_{q,p}(sl(N))$ exchange algebras}

\subsection{$q$-Virasoro exchange algebras}

While for $c=-N$, $t(z)$ and $L(w)$ commute, for 
$(-p^{\frac{1}{2}})^{NM} = q^{-c-N}$ $(M \in {\mathbb Z})$, $t(z)$ and 
$L(w)$ realize an \emph{exchange algebra}:
\[
t(z) \, L(w) = F(M,w/z) \, L(w) \, t(z)
\]
where
\[
F(M,x) = \left\{ \begin{array}{ll} 
\displaystyle q^{2M(N-1)} \prod_{k=0}^{NM-1}
\frac{\Theta_{q^{2N}}(x^{-2}p^{-k}) \, \Theta_{q^{2N}}(x^2p^{k})}
{\Theta_{q^{2N}}(x^{-2}q^2p^{-k}) \, \Theta_{q^{2N}}(x^2q^2p^{k})}
& \mbox{for $M>0$} \,, \\ \\
\displaystyle q^{-2|M|(N-1)} \prod_{k=1}^{N|M|}
\frac{\Theta_{q^{2N}}(x^{-2}q^2p^{k}) \, \Theta_{q^{2N}}(x^2q^2p^{-k})}
{\Theta_{q^{2N}}(x^{-2}p^{k}) \, \Theta_{q^{2N}}(x^2p^{-k})}
& \mbox{for $M<0$} \,. \\
\end{array} \right.
\]
It follows that one gets the following theorem: \\
\textbf{Theorem 3:}  
On the 2d-surface $\Sigma_{N,M} = \{ (p,q) \, \vert \, 
(-p^{\frac{1}{2}})^{NM} = q^{-c-N} \}$ of the 3d-space of parameters 
$q,p,c$, the operators generated by $t(z)$ close a quadratic algebra:
\[
t(z) \, t(w) = {\cal Y}_{N,p,q,M}\Big(\frac{w}{z}\Big) \, t(w) \, t(z)
\]
where
\[
{\cal Y}_{N,p,q,M}(x) = \left\{ \begin{array}{ll}
\displaystyle \prod_{k=1}^{NM} \frac{\Theta_{q^{2N}}^2(x^2 p^{-k}) \,
\Theta_{q^{2N}}(x^2 q^2 p^k) \, \Theta_{q^{2N}}(x^2 q^{-2} p^k)}
{\Theta_{q^{2N}}^2(x^2 p^k) \, \Theta_{q^{2N}}(x^2 q^2 p^{-k})
\, \Theta_{q^{2N}}(x^2 q^{-2} p^{-k})} & \mbox{for $M>0$} \,, \\ \\
\displaystyle \prod_{k=1}^{N|M|-1}
\frac{\Theta_{q^{2N}}^2(x^2 p^{-k}) \,
\Theta_{q^{2N}}(x^2 q^2 p^k) \, \Theta_{q^{2N}}(x^2 q^{-2} p^k)}
{\Theta_{q^{2N}}^2(x^2 p^k) \, \Theta_{q^{2N}}(x^2 q^2 p^{-k})
\, \Theta_{q^{2N}}(x^2 q^{-2} p^{-k})} & \mbox{for $M<0$} \,. \\
\end{array} \right.
\]
{For} $N = 2$, one has
\begin{eqnarray*}
&& {\cal Y}_{2,p,q,M}(xq^2) = {\cal Y}_{2,p,q,M}(x) \,, \\
&& {\cal Y}_{2,p,q,M}(x) \, {\cal Y}_{2,p,q,M}(xq) = 1 \,.
\end{eqnarray*}
One recovers the same periodicity properties satisfied by structure 
function $\varphi(x)$ of \cite{SKAO} and \cite{FF} in the case of 
$sl(2)$ (when $M=1$, one makes the change of variable $q^2 \rightarrow 
p$, $p \rightarrow q$, $x^2 \rightarrow x$, then ${\cal Y}_{2,p,q,1}(x) 
\rightarrow \varphi^2(x)$ of \cite{FF}).  \\
{For} $N > 2$, one has
\begin{eqnarray*}
&& {\cal Y}_{N,p,q,M}(xq^N) = {\cal Y}_{N,p,q,M}(x) \,, \\
&& {\cal Y}_{N,p,q,M}(x) \, {\cal Y}_{N,p,q,M}(xq) \, \dots 
\, {\cal Y}_{N,p,q,M}(xq^{N-1}) = 1 \,.
\end{eqnarray*}

\subsection{${\cal W}_{q,p}(sl(N))$ exchange algebras}

The generalization to higher spin generators in the quantum case follows 
that of the classical case, although with a required notion of ordering 
between individual $t(z)$-generators:
\[
s_i(z) = \prod_{\frac{i-1}{2} \ge u \ge -\frac{i-1}{2}}^{{\textstyle 
\curvearrowleft}} t(q^uz) \,.
\]
The exchange algebra takes now the form:
\[
s_i(z) s_j(w) = \prod_{u=-\frac{i-1}{2}}^{\frac{i-1}{2}} 
\prod_{v=-\frac{j-1}{2}}^{\frac{j-1}{2}}{\cal Y}_{N,p,q,M} 
\left(q^{v-u}\frac{w}{z}\right) \, s_j(w) \, s_i(z) \,.
\]
Furthermore one can verify that $s_N(z)$ commutes with all other 
generators; hence we are justified in restricting $i$ to $1, \ldots , 
N-1$.

\subsection{Classical limit}

\textbf{Theorem 4:} 
On the surface $\Sigma_{N,M}$, when $p = q^{Nh}$ with $h \in {\mathbb Z} 
\backslash \{ 0 \}$:
\[
{\cal Y}_{N,p,q,M} = 1 \,.
\]
Hence $t(z)$ realizes an Abelian subalgebra in ${\cal 
A}_{q,p}(\widehat{sl}(N)_{c})$ (though it is not in the center of ${\cal 
A}_{q,p}(\widehat{sl}(N)_{c})$ in general).

The result of the previous theorem allows us to define Poisson 
structures on the corresponding Abelian algebras generated by $t(z)$.  
They are obtained as limits of the exchange algebra of Theorem 3 when $p 
= q^{Nh}$ with $h \in {\mathbb Z} \backslash \{ 0 \}$.  Setting $q^{Nh} = 
p^{1-\beta}$, the $h$-labeled Poisson structure is defined by:
\[
\{ t(z) , t(w) \}^{(h)} \equiv \lim_{\beta \rightarrow 0} \frac{1}{\beta}
\, ( t(z)t(w) - t(w)t(z) ) \,.
\]
{For} $h$ even, $\{ t(z) , t(w) \}^{(h)}$ coincides with the Poisson 
structure of the critical case $c=-N$, while for $h$ odd, $\{ t(z) , 
t(w) \}^{(h)}$ gives new Poisson structures (for $sl(N)$ when $N > 2$).
Explicitly, one has the following expression:
\[
\{ t(z) , t(w) \}^{(h)} = f_h(w/z) \, t(z) \, t(w)
\]
with
\begin{eqnarray*}
f_h(x) &=&
{\cal N}_{odd} \left[ \sum_{\ell \ge 0}
E(\sfrac{NM}{2})(E(\sfrac{NM}{2})+1)
\left( \frac{2x^2q^{2N\ell}}{1-x^2q^{2N\ell}}
- \frac{x^2q^{2N\ell+2}}{1-x^2q^{2N\ell+2}} \right. \right. \\
&& \left. - \frac{x^2q^{2N\ell-2}}{1-x^2q^{2N\ell-2}} \right)
+ E(\sfrac{NM+1}{2})^2 \left(
\frac{2x^2q^{2N\ell+N}}{1-x^2q^{2N\ell+N}}
- \frac{x^2q^{2N\ell+N+2}}{1-x^2q^{2N\ell+N+2}} \right. \\
&& \left. - \frac{x^2q^{2N\ell+N-2}}{1-x^2q^{2N\ell+N-2}} \right)
- \sfrac{1}{2} E(\sfrac{NM}{2})(E(\sfrac{NM}{2})+1)
\left( \frac{2x^2}{1-x^2} - \frac{x^2q^2}{1-x^2q^2} \right. \\
&& \left.  \left.  - \frac{x^2q^{-2}}{1-x^2q^{-2}} \right) - (x 
\leftrightarrow x^{-1}) \right] \qquad \qquad \qquad \qquad \qquad 
\qquad \quad \mbox{for $h$ odd} \\
&& \\
&=& {\cal N}_{even} \left[ \sum_{\ell \ge 0} \left(
\frac{2x^2q^{2N\ell}}{1-x^2q^{2N\ell}}
- \frac{x^2q^{2N\ell+2}}{1-x^2q^{2N\ell+2}}
- \frac{x^2q^{2N\ell-2}}{1-x^2q^{2N\ell-2}} \right) \right. \\
&& \left. - \sfrac{1}{2} \left( \frac{2x^2}{1-x^2}
- \frac{x^2q^2}{1-x^2q^2} - \frac{x^2q^{-2}}{1-x^2q^{-2}} \right)
- (x \leftrightarrow x^{-1}) \right]
\qquad \mbox{for $h$ even} \\
\end{eqnarray*}
where ${\cal N}_{odd} = 2Nh \ln q$, ${\cal N}_{even} = N^2M(NM+1)h \ln 
q$ and the notation $E(n)$ means the integer part of the number $n$.

\medskip

The mode Poisson brackets in the sector $k=0$ is given by:
\begin{eqnarray*}
\{ t_n , t_m \}^{(h)} &=&
-2(q-q^{-1}) \ln q \sum_{r \in {\mathbb Z}} \left(
- E(\sfrac{NM+1}{2})^2 \frac{[r]_{q}^2}{[Nr]_{q}} \right. \\
&& \left. + E(\sfrac{NM}{2})(E(\sfrac{NM}{2})+1)
\frac{[(N-1)r]_{q}[r]_{q}}{[Nr]_{q}} \right) \, t_{n-2r} \, t_{m+2r} \\
&& \hspace{70mm} \mbox{for $h$ odd} \\
&& \nonumber \\
&=& -2(q-q^{-1}) \ln q \sum_{r \in {\mathbb Z}}
\frac{[(N-1)r]_{q}[r]_{q}}{[Nr]_{q}} \, t_{n-2r} \, t_{m+2r} \\
&& \hspace{70mm} \mbox{for $h$ even}
\end{eqnarray*}
Since the initial non-abelian structure for $t(z)$ is closed, the 
exchange algebras
\[
t(z) \, t(w) = {\cal Y}_{N,p,q,M}(w/z) \, t(w) \, t(z)
\]
are \emph{inequivalent quantizations} with respect to $p$ of the Poisson 
bracket which defines $Vir_q(sl(N))$: $t(z)$ generates a 
$Vir_{q,p}(sl(N))$ algebra.

\medskip

The generalization for the higher spin generators $s_{i}(z)$ presents no 
conceptual difficulty. The following results are obtained: \\
-- for $h$ even, $\{ s_{i}(z) , s_{j}(w) \}^{(h)}$ coincides with the 
Poisson structure of the critical case $c=-N$.  The Poisson bracket is 
again identical to the ``core'' contribution in \cite{FR} (i.e.  without 
lower-spin extensions).  \\
-- for $h$ odd, $\{ s_{i}(z) , s_{j}(w) \}^{(h)}$ gives rise to two new 
types of quadratic Poisson structures (that is $q$-deformed classical 
${\cal W}_{N}$-algebras) whether $i+j \le N$ or $i+j > N$.  \\

\section{Conclusion}

Here above is summarized the construction, from elliptic algebras, of 
the deformed Virasoro and ${\cal W}$ algebras in the classical case as 
well as in the quantum case.  We can mention that a study for consistent 
central extensions, providing in particular general results for the 
simplest cases of $Vir_{q}(sl(2))$ has very recently been achieved in 
ref.  \cite{AFRS4}.  Also, a direct determination of abstract generators 
for $q$-deformed ${\cal W}_{N}$ algebras has been carried out.  Such an 
universal construction of ${\cal W}_{q,p}$ algebras is achieved in 
\cite{AFRS5}.

\medskip

An important application of the above developed constructions is the 
determination of integrable models submitted to such symmetries.  It is 
known that the trigonometric Calogero--Moser model is related to the 
Virasoro algebra and the Ruijsenaars--Schneider model to the 
$q$-deformed Virasoro algebra considered in ref.  \cite{SKAO}.  One 
expects the existence of generalized Ruijsenaars--Schneider models the 
symmetries of which would be given by the just described $q$-Virasoro 
exchange algebras.  For such a program, the bosonization of these 
algebras appears as an essential tool, which is still missing.

\end{document}